\documentclass[11pt]{amsart}
\usepackage{amsmath, amsthm, mathabx,amscd, amsfonts, amssymb, hyperref, mathrsfs, textcomp, epsfig, a4wide, graphicx, IEEEtrantools, tikz, verbatim, xypic }

\usepackage[all,cmtip]{xy}
\usetikzlibrary{matrix, decorations.pathreplacing,angles,quotes,cd}

\newtheorem{thm}{Theorem}[section]

\newtheorem{prop}[thm]{Proposition}

\newtheorem{conj}[thm]{Conjecture}
\theoremstyle{definition}
\newtheorem{defn}[thm]{Definition}

\usepackage{mathtools}

\DeclarePairedDelimiter\floor{\lfloor}{\rfloor}

\newcommand{\spin}{\mathfrak{s}}

     \RequirePackage{rotating}                   
    \def\HMt{%
       \setbox0=\hbox{$\widehat{\mathit{HM}}$}
       \setbox1=\hbox{$\mathit{HM}$}
       \dimen0=1.1\ht0
       \advance\dimen0 by 1.17\ht1
       \smash{\mskip2mu\raise\dimen0\rlap{%
          \begin{turn}{180}
              {$\widehat{\phantom{\mathit{HM}}}$}
           \end{turn}} \mskip-2mu    
                \mathit{HM}
    }{\vphantom{\widehat{\mathit{HM}}}}{}}

\newcommand{\HMb}{\overline{\mathit{HM}}}
\newcommand{\HMf}{\widehat{\mathit{HM}}}

\theoremstyle{remark}
\newtheorem{remark}{Remark}[section]

\begin{document}
\title[Families of Dirac operators]{Lectures on families of Dirac operators and applications}

\author{Francesco Lin}
\address{Department of Mathematics, Columbia University} 
\email{flin@math.columbia.edu}

\maketitle

 \vspace{0.5cm}

 {\centering\textit{Dedicated to Tom Mrowka in the occasion of his 60th birthday.}\par}

\begin{abstract}These are the notes for a minicourse taught at the 2022 ICTP summer school \textit{Frontiers in Geometry and Topology}. The goal is to introduce families of Dirac operators and how they can be used to study interactions between geometry and topology. In particular, we discuss the index theorem for families of twisted Dirac operators, and its applications to metrics of positive scalar curvature and the three-dimensional Weinstein conjecture.

\end{abstract}

\section*{Introduction}

These are the notes for the three and a half hour long minicourse \textit{Families of Dirac operators and applications} taught at the 2022 ICTP summer school \textit{Frontiers in Geometry and Topology} for an audience of graduate students. The main goal is to show via concrete examples how families of Dirac operators can be used as a bridge between the geometry and the topology of manifolds. 
These ideas are closely related to current topics of research in gauge theory and low-dimensional topology, and the hope is to provide a friendly first introduction to the subject for those working in the field. In particular, these notes assume the reader to be somehow familiar with the following:
\begin{itemize}
\item vector bundles, their classifying spaces, and characteristic classes \cite{MilSta};
\item basic differential geometry, including connections on vector bundles and the various notions of curvature of a Riemannian manifold \cite{Lee};
\item basic notions in functional analysis such as separable Hilbert spaces and bounded operators \cite[Appendix D]{Eva}.
\end{itemize}
The last part of the notes will also assume the reader to be familiar with the formal properties of either monopole or Heegaard Floer homology (see for example \cite{Lin} and \cite{Gre} respectively for introductions to the topics).
\\
\par
The minicourse consisted of three lectures covering the following:
\begin{enumerate}
\item an informal introduction to spin Dirac operators; the Atiyah-Singer index theorem, the Lichnerowicz formula, and applications to metrics of positive scalar curvature.
\item $K$-theory and its relation with Fredholm operators; the index theorem for families of Dirac operators and the Gromov-Lawson theorem on the torus.
\item the space of self-adjoint Fredholm operators; relations with Taubes' proof of the Weinstein conjecture in dimension three.
\end{enumerate}
This is a very ambitious syllabus - each lecture could have been the topic of a minicourse itself. When preparing the lectures, we have tried to highlight the motivation and key ideas behind the constructions while keeping technicalities at the bare minimum (some important ones are pointed out in the footnotes); we provide precise references for more detailed treatments. These notes follow the content of the lectures very closely, and will have a rather terse style. At the end of each lecture, the interested reader will find some exercises.

\subsection*{Warning.} Throughout the text, we will be very imprecise about the analytic setup. For example, we will state most index theorems in terms of bounded operators, but switch to unbounded (differential) ones without further notice when discussing geometric applications.

\subsection*{Acknowledgements.} The author thanks the participants of the summer school, whose questions greatly helped to improve these notes. This work was partially supported by the Alfred P. Sloan Foundation and NSF grant DMS-2203498.

\vspace{2cm}

\section{Dirac operators, geometry and topology}

\subsection{Dirac operators on $\mathbb{R}^n$}
The most fundamental differential operator on $\mathbb{R}^n$ is the Laplacian\footnote{Our convention is to use all negative signs so that $\langle u,\Delta u\rangle_{L^2}\geq 0 \text{ for all compactly supported $u$}$, where we consider the $L^2$-inner product $\langle u,v\rangle_{L^2}=\int uv d\mathrm{vol}$.}
\begin{equation*}
\Delta:=-\frac{\partial^2}{\partial x_1^2}-\frac{\partial^2}{\partial x_2^2}-\dots-\frac{\partial^2}{\partial x_n^2}.
\end{equation*}
Functions satisfying $\Delta u=0$ are called \textit{harmonic}; these can be used to describe many important objects, such as potentials in electromagnetism. While working on the relativistic quantum mechanics of the electron, Dirac \cite{Dir} was led to look for a first order differential operator $D$ such that $D^2=\Delta$.
\begin{remark}
To be precise, Dirac was working with a Lorenzian metric on $\mathbb{R}^4$, so that he considered the d'Alembertian 
\begin{equation*}
\Box=-\frac{\partial^2}{\partial x_1^2}-\frac{\partial^2}{\partial x_2^2}-\frac{\partial^2}{\partial x_3^2}+\frac{\partial^2}{\partial t^2}
\end{equation*}
instead of the Laplacian $\Delta$. We will always work with Riemannian metrics in these notes.
\end{remark}
Heuristically, what does such an operator $D$ look like? Let us assume for simplicity that it takes the form
\begin{equation*}
D=a_1 \frac{\partial}{\partial x_1}+a_2 \frac{\partial}{\partial x_2}+\dots +a_n \frac{\partial}{\partial x_n}
\end{equation*}
with the $\{a_i\}$ are constant quantities (belonging to some algebraic structure to be decided). Then
\begin{equation*}
D^2=\sum_i a_i^2\frac{\partial^2}{\partial x_i^2}+\sum_{i,j} a_ia_j\frac{\partial}{\partial x_i}\frac{\partial}{\partial x_j}.
\end{equation*}
Hence, as derivatives commute, $D^2=\Delta$ is equivalent to the relations
\begin{equation}\label{clifford}
\begin{cases}
a_i^2&=-1\\
a_ia_j+a_ia_j&=0 \text{ if }i\neq j.
\end{cases}
\end{equation}
Let us see how to solve such system of equations in the $\{a_i\}$. If $n=1$ we just need $a_1^2=-1$; of course, we cannot solve this within real numbers, but if we allow complex valued functions, we can set $a_1=i$, and
\begin{equation*}
D=i\frac{\partial}{\partial x_1}
\end{equation*}
squares to $\Delta$ (acting now on complex-valued functions).
\par
When $n=2$, the anticommutativity relation $a_1a_2+a_2a_1=0$ implies that we cannot work anymore with complex numbers; on the other hand, we can solve (\ref{clifford}) if we allow the $a_i$ to be $2\times 2$ matrices by setting
\begin{equation*}
a_1=\begin{bmatrix}
0&-1\\
1&0
\end{bmatrix}
\quad
a_2=\begin{bmatrix}
0&i\\
i&0
\end{bmatrix},
\end{equation*}
so that
\begin{equation}\label{DiracR2}
D=a_1 \frac{\partial}{\partial x_1}+a_2 \frac{\partial}{\partial x_2}=\begin{bmatrix}
0&-\frac{\partial}{\partial x_1}+i \frac{\partial}{\partial x_2}\\
\frac{\partial}{\partial x_1}+i \frac{\partial}{\partial x_2}&0
\end{bmatrix}
\end{equation}
squares to $\Delta$ thought of as an operator on $\mathbb{C}^2$ valued functions. Notice that $D$ is closely related to the holomorphic and anti-holomorphic derivative operators.
\par
In higher dimensions, the situation is similar; when $n=3$, we can set the $a_i$'s to be the Pauli matrices
\begin{equation*}
\sigma_1=
\begin{bmatrix}
i&0\\
0&-i
\end{bmatrix}
\quad
\sigma_2=
\begin{bmatrix}
0&-1\\
1&0
\end{bmatrix}
\quad
\sigma_3=
\begin{bmatrix}
0&i\\
i&0
\end{bmatrix},
\end{equation*}
and the corresponding Dirac operator $D$ is again a square root of $\Delta$ acting on $\mathbb{C}^2$ valued functions. When $n=4$, we can set the $a_i$ to be the $4\times 4$ block matrices defined by
\begin{equation*}
\begin{bmatrix}
0&-I_2\\
I_2&0
\end{bmatrix}
\text{ and }
\begin{bmatrix}
0&-\sigma_i^*\\
\sigma_i&0
\end{bmatrix}\text{ for }i=1,2,3
\end{equation*}
where the $I_2$ is the $2\times 2$ identity matrix and $^*$ denotes transpose conjugate. The corresponding Dirac operator $D$ squares to $\Delta$ acting on $\mathbb{C}^4$ valued functions, and is the (Riemannian version of the) operator Dirac was looking for. We will leave to the reader to guess the general recipe to define $D$ on $\mathbb{R}^n$ for $n\geq 5$.
\begin{remark}
Of course, we could have made different choices for the matrices $a_i$, which would have lead to different Dirac operators; for example, when $n=1$ we could have taken $D=i\frac{d}{dx_1}$ to be acting on $\mathbb{C}^m$-valued functions. The specific operators we constructed above are `minimal' in a specific sense (involving the representation theory of \textit{Clifford algebras}), and are called the \textit{spin Dirac operators}.
\end{remark}

\vspace{0.3cm}

\subsection{Dirac operators on a manifold}
In this lectures we will be interested in how spin Dirac operators can be used as a bridge between the geometry and topology of an oriented Riemannian manifold $(M,g)$; the precise definition of the analogue $D$ is quite involved, and we refer the reader for example to \cite[Ch. 3,4]{Roe} for more details. We will only focus on formal properties of such generalization.
\par
First of all, in order to define $D$, we need to choose a \textit{spin structure} on $M$; we will not need the precise definition of what this is (see for example \cite[Ch. IV]{Kir}), but will only remark that this exists if and only the second Stiefel-Whitney class $w_2(TM)$ is zero. For manifolds of dimension $n\geq3$, this is equivalent to the restriction of $TM$ to the $2$-skeleton of $M$ being trivial. A spin structure affords us an analogue of the target vector spaces $\mathbb{C}^m$ we considered before, which is a hermitian bundle
\begin{equation*}
S\rightarrow M
\end{equation*}
of rank $2^{\floor{n/2}}$ equipped with a hermitian connection $\nabla$, called the spinor bundle. The spin Dirac operator is then a first order differential operator
\begin{equation*}
D:\Gamma(S)\rightarrow \Gamma(S)
\end{equation*}
with the following key properties:
\begin{itemize}
\item it is \textit{formally self-adjoint in $L^2$}, i.e. if $\varphi,\psi\in\Gamma(S)$ are compactly supported sections (i.e. if $M$ is compact), then
\begin{equation*}
\langle D\varphi,\psi\rangle_{L^2}=\langle \varphi,D\psi\rangle_{L^2}.
\end{equation*}
\item it is \textit{elliptic}; we will not define exactly what this means, but this implies for example that \textit{harmonic spinors}, i.e. solutions to the equation $D\varphi=0$ (which makes sense for $\varphi$ merely $C^1$) are automatically smooth. \end{itemize}
One should compare the latter with the analogous result for holomorphic and harmonic functions, which are in the kernel of the elliptic operators $\overline{\partial}$ and $\Delta$ respectively. We refer the reader to \cite[Chapter 6]{War} for a nice introduction to the analysis of elliptic equations on manifolds which only assumes some basic familiarity with Fourier series.

\vspace{0.3cm}

\subsection{Dirac operators and topology} Another key consequence of ellipticity is that, when $M$ is compact, the space of harmonic spinors is \textit{finite} dimensional. This is surprisingly related to topological information about $M$ in a way that we now discuss.
\par
Assume that $n$ is even. Looking back at the case of $\mathbb{R}^n$, we see that in this case the matrices corresponding to $a_i$ have a simple structure as block matrices. In the case of a Riemannian manifold, this manifests in a canonical splitting of the spinor bundle
\begin{equation*}
S=S^+\oplus S^-
\end{equation*}
for which we can write the Dirac operator in block form
\begin{equation*}
D=\begin{bmatrix}
0&D^-\\
D^+&0
\end{bmatrix}
\text{ where }
D^{\pm}:\Gamma(S^{\pm})\rightarrow \Gamma(S^{\mp}).
\end{equation*}
We call $D^{\pm}$ the \textit{chiral} Dirac operators. As $D$ is self-adjoint, $D^{\pm}$ are adjoint (i.e. transpose) to each other; in particular, the analogue in this situation of the relation $(\mathrm{im}A)^{\perp}=\mathrm{ker}A^T$ for a matrix $A$ implies that
\begin{equation*}
\mathrm{dim}\mathrm{coker} D^+=\mathrm{dim}\mathrm{ker} D^-,
\end{equation*}
hence is also finite dimensional\footnote{There is a subtlety here: we use that $\mathrm{im}D^+$ is a closed subspace to identify its orthogonal with $\mathrm{coker} D^+$; closedness also follows from ellipticity.}. In general, an operator $T$ between Hilbert spaces with finite dimensional kernel and cokernel is called \textit{Fredholm}; its \textit{index} is defined to be
\begin{equation*}
\mathrm{ind}(T)=\mathrm{dim}\mathrm{ker}T-\mathrm{dim}\mathrm{coker}T\in\mathbb{Z}.
\end{equation*}
The following fundamental result (a special case of the Atiyah-Singer index theorem, one of the cornerstones of 20th century mathematics) relates the index of $D^+$ (which is an analytical quantity) to the topology of $M$; this is particularly striking because $\mathrm{dim}\mathrm{ker} D^+$ and $\mathrm{dim}\mathrm{coker} D^+$ are not themselves topological invariants (cf. \cite{Hit})!
\begin{thm}[\cite{AS}]\label{AS}
The index of $D^+$ can be computed in terms of the Pontryagin classes $p_i(TM)$. In particular:
\begin{itemize}
\item if $n\equiv 2$ modulo $4$, $\mathrm{ind}(D^+)=0$;
\item if $n=4$,  $\mathrm{ind}(D^+)=-\langle p_1,[M]\rangle/24$;
\item if $n=8$, $\mathrm{ind}(D^+)=\langle -4p_2+7p_1^2,[M]\rangle/5760$.
\end{itemize}
More in general, the index is determined by an expression in the Pontryagin classes of $TM$ called $A$-hat genus of $M$, and $\mathrm{ind}(D^+)=\langle \widehat{\mathcal{A}}(M),[M]\rangle$.
\end{thm}
\begin{remark}When $n=4$, the Hirzebruch signature formula $p_1=3\sigma$ implies that we can also write $\mathrm{ind}(D^+)=-\sigma(M)/8$, where $\sigma$ denotes the signature.
\end{remark}
For a more thorough discussion of the Atiyah-Singer index theorem and its various proofs we refer the reader to \cite{Sha} and \cite{Roe}.

\vspace{0.3cm}

\subsection{Dirac operators and geometry} Recall that Dirac's goal was to find a differential operator such that $D^2=\Delta$. The Dirac operator on a spin manifold satisfies an analogous relation. Indeed, recall that our bundle $S\rightarrow M$ comes with a hermitian connection
\begin{equation*}
\nabla:\Gamma(S)\rightarrow \Gamma(T^*M\otimes S).
\end{equation*}
We denote the formal adjoint in $L^2$ of such operator by
\begin{equation*}
\nabla^*: \Gamma(T^*M\otimes S)\rightarrow\Gamma(S).
\end{equation*}
This is the differential operator determined by the identity
\begin{equation*}
\langle \nabla\varphi,\alpha\rangle_{L^2}=\langle \varphi,\nabla^*\alpha\rangle_{L^2}\text{ for every }\varphi\in\Gamma(S),\alpha\in \Gamma(T^*M\otimes S),
\end{equation*}
and is (after changing sign) the bundle analogue of the classical divergence of a vector field.
\par
The second order operator $\nabla^*\nabla$ is called the \textit{Bochner Laplacian}; it makes sense on any bundle with Hermitian connection, and corresponds in $\mathbb{R}^n$ to the familiar expression
\begin{equation*}
\Delta=-\mathrm{div}\circ\mathrm{grad}.
\end{equation*}
The fundamental identity relating the Dirac operator and geometry is then the following.
\begin{thm}[Lichnerowicz `63]\label{Lichnerowicz} Let $s$ be the scalar curvature of $(M,g)$. Then $D^2=\nabla^*\nabla+s/4$.
\end{thm}
The magic behind the formula is the following: $D^2$ and $\nabla^*\nabla$ are second order differential operators which by design have the same second order part; one would expect then their difference to be a first order differential operator, but it turns out to be a zeroth order one! This has surprising consequences for the study of the existence of metrics of positive scalar curvature, which is a quite subtle problem in Riemannian geometry.

\begin{prop}\label{Bochner}Suppose $(M,g)$ is a compact spin manifold with $s>0$ everywhere. Then the identity $\mathrm{ind}(D^+)=0$ (which is topological in nature by Theorem \ref{AS}) holds. In particular, if $M$ is $4$-dimensional, $\sigma(M)=0$.
\end{prop}
So for example the $K3$ surface, which is a spin $4$-manifold with $\sigma=-16$ (see \cite{GS}), does not admit metrics of positive scalar curvature. Notice that the spin assumption is essential, as the Fubini-Study metric on $\mathbb{C}P^2$ has $s>0$.

\par
To prove the proposition, we just need to show that $D$ (and therefore $D^{\pm}$) has trivial kernel. To see this, assume $D\varphi=0$ and $\varphi$ is not identically zero. The Lichnerowicz formula then implies
\begin{equation*}
\nabla^*\nabla\varphi+\frac{s}{4}\varphi=0.
\end{equation*}
Taking the $L^2$-inner product with $\varphi$, we obtain
\begin{equation*}
0=\langle\nabla^*\nabla\varphi,\varphi\rangle_{L^2}+\langle \frac{s}{4}\varphi,\varphi\rangle_{L^2}>\langle\nabla^*\nabla\varphi,\varphi\rangle_{L^2}=\langle\nabla\varphi,\nabla\varphi\rangle_{L^2}\geq0
\end{equation*}
which is a contradiction.

\vspace{0.3cm}

\subsection{Twisted Dirac operators}

Consider now the case of the torus $T^n$, which admits a flat metric, hence a metric with $s\equiv 0$; it is a natural ask whether it admits a metric with $s>0$. The answer is negative for $n=2$ by Gauss Bonnet's theorem. Also, notice that Myers' theorem \cite[Chapter $12$]{Lee} implies that $T^n$ does not admit metrics with positive Ricci curvature because it has infinite fundamental group.
\par
 The tangent bundle of ${T}^n$ is trivial, so in particular it is spin ($w_2=0$) and has vanishing Pontryagin classes; hence $\mathrm{ind}(D^+)=0$ and the obstruction of Proposition \ref{Bochner} does not apply.
\par
The key idea to obtain more subtle obstructions is to study twisted Dirac operators. In particular given $M$, consider flat hermitian connections $B$ on the trivial line bundle $\underline{\mathbb{C}}\rightarrow M$; up to isomorphism, these are classified by their holonomy around loops, and correspond bijectively to the $b_1(M)$-dimensional torus
\begin{equation*}
\mathbb{T}_M=H^1(M;\mathbb{R})/H^1(M;2\pi \mathbb{Z}).
\end{equation*}
To each of these, one can associate a twisted connection $\nabla_{B}$ on $S$, and correspondingly a twisted Dirac operator $D_B$ on $\Gamma(S)$ which satisfies the analogues of Theorem \ref{AS} and \ref{Lichnerowicz}\footnote{Here it is important that $B$ is flat.}. For example, on $M=S^1=\mathbb{R}/2\pi \mathbb{Z}$, the twisted Dirac operators are of the form
\begin{equation*}
i\frac{d}{dt}+c \text{ for }c\in\mathbb{R}.
\end{equation*}
Two such operators are conjugate exactly when $c-c'\in\mathbb{Z}$, so that we get a $\mathbb{T}_{S^1}=S^1$-family of operators. In the next lecture, we will see that the family of Fredholm operators $\{D_B\}$ carries much more information than a single operator.

\vspace{0.3cm}

\subsection{Exercises}
\begin{enumerate}
\item The Dirac operator (\ref{DiracR2}) on $\mathbb{R}^2$ is translation invariant, and descends to the Dirac operator on the torus $\mathbb{R}^2/\mathbb{Z}^2$. Determine the space of harmonic spinors, and check that the index formula for $D^+$ holds in this case. 
\item By looking carefully at the argument in Proposition $1.3$, show that if $(M,g)$ is a spin $4$-manifold with $s\equiv0$, then $|\sigma(M)|\leq 16$.
\item Use Fourier series to compute the spectrum of the twisted Dirac operator $i\frac{d}{dt}+c$ on $S^1$; show that two such operators are conjugate when $c-c'\in\mathbb{Z}$.
\end{enumerate}

\vspace{2cm}

\section{K-theory and the index theorem for families}

Consider a complex separable Hilbert space $H$, e.g. $H=L^2(M)$. Recall that a bounded operator $T:H\rightarrow H$ is Fredholm if it has finite dimensional kernel and cokernel\footnote{Because we are working with Hilbert (hence Banach) spaces, the image of such $T$ is automatically closed; when working with more general topological vector spaces the latter needs to be incorporated in the definition of Fredholm}. Denote by $\mathrm{Fred}(H)$ the space of such operators; one can show that it is an open set of the space of bounded operators, and the map
\begin{equation*}
\mathrm{ind}:\mathrm{Fred}(H)\rightarrow \mathbb{Z}
\end{equation*}
is locally constant (even though the dimensions of kernel and cokernel are not!). It turns out that the space $\mathrm{Fred}(H)$ has very interesting topology, which we now describe.

\subsection{K-theory}
We discuss the essential notions in topological $K$-theory, and refer the reader to \cite[Ch. 38]{FF} for more details. Consider a finite CW complex $X$, which will assume for simplicity to be connected (in our applications $X$ will be the torus $\mathbb{T}_M$). Denote by $\mathrm{Vect}(X)$ the set of isomorphism classes of complex vector bundles on $X$. The operation of direct sum $\oplus$ and tensor product $\otimes$ make $\mathrm{Vect}(X)$ into a \textit{semiring}, i.e. a structure like a ring but without additive inverses. For example, $\mathbb{N}$ is a semiring, and there is indeed a semiring homomorphism $\mathrm{Vect}(X)\rightarrow \mathbb{N}$ given by the rank\footnote{If $X$ is not connected, one allows vector bundles with different ranks on each of them, and this map depends on a choice of connected component.}. Notice that, unlike $\mathbb{N}$, $\mathrm{Vect}(X)$ is not in general \textit{cancellative}, i.e. $a+c=b+c$ does not in general imply $a=b$.
\par
Of course, it is nicer in general to work with rings than semirings. Given any semiring $S$, there is a natural way to associate a ring called the \textit{Grothendieck ring} $\mathrm{Gr}(S)$. This follows the usual construction of $\mathbb{Z}$ from $\mathbb{N}$; one defines $\mathrm{Gr}(S)$ to be the set of formal differences $a-b$ of elements in $S$, modulo the relation
\begin{equation*}
a-b\sim a'-b'\text{ if }a+b'+e=a'+b+e\text{ for some }e.
\end{equation*}
Here we need the extra $e$ to have a good definition for non-cancellative semirings.
\begin{defn}
For a finite CW complex $X$, the $K$-theory of $X$, denoted by $K(X)$, is the Grothendieck ring of $\mathrm{Vect}(X)$.
\end{defn}
Here are some basic properties:
\begin{itemize}
\item An element of $K(X)$ is represented by a formal difference of vector bundles $E-F$;
\item There is a ring homomorphism $K(X)\rightarrow \mathbb{Z}$ given by the rank;
\item $K$ gives rise to a contravariant functor from finite CW complexes to rings; homotopic maps induce the same ring homomorphism.
\end{itemize}
Denote by $\mathrm{Vect}_k(X)$ the isomorphism classes of $k$-dimensional vector bundles. Recall that we have the identification with the set of homotopy classes
\begin{equation*}
\mathrm{Vect}_k(X)=[X,BU(k)]
\end{equation*}
where, concretely, $BU(k)$ is the union of the Grassmannians $G(n,k)$ of $k$-planes in $\mathbb{C}^n$. Now, for any vector bundle $E\rightarrow X$, there is a vector bundle $E'$ such that $E\oplus E'$ is trivial. This implies that any element in $K(X)$ has a representative of the form $F-\underline{\mathbb{C}}^m$ for some vector bundle $F$ and $m\in\mathbb{N}$. From this, one shows that
\begin{equation*}
K(X)=[X,BU\times\mathbb{Z}]
\end{equation*}
where $BU=\bigcup_{k\geq1} BU(k)$, and $\mathbb{Z}$ corresponds to the rank.

\vspace{0.3cm}

\subsection{Relation with Fredholm operators}

It turns out that $K$-theory is the right algebraic tool to study families of Fredholm operators because of the following.
\begin{thm}[Atiyah-J\"anich] For a finite CW complex $X$ , there is a natural bijection
\begin{equation*}
K(X)=[X,\mathrm{Fred}(H)],
\end{equation*}
i.e. the space of Fredholm operators is also a classifying space for $K$-theory. In particular, homotopy classes of families of Fredholm operators parametrized by $X$ are in bijection with $K(X)$.
\end{thm}
As a consequence, this also tells us that the space $\mathrm{Fred}(H)$ has a very interesting topology.
\par
Following \cite[Appendix]{Ati}, let us discuss how to associate to a family $\{T_x\}_{x\in X}$ of Fredholm operators $T_x:H\rightarrow H$ a class in $K(X)$. Intuitively, this corresponds to the difference between the kernel and cokernel at each point of $T_x$. Of course, these do not form a vector bundle (the dimensions might not even be constant!) so we need a different approach to formalize this. Because $X$ is compact, one can find a finite dimensional subspace $W\subset H$ such that
\begin{equation*}
\mathrm{Im}(T_x)+W=H
\end{equation*}
for all $x\in X$. Then, $V=\bigcup T_x^{-1}(W)$ is a vector bundle over $X$, and we associate to $\{T_x\}$ the element
\begin{equation*}
\mathrm{ind}\{T_x\}:=V-\underline{\mathbb{C}}^{\mathrm{dim}W}\in K(X)
\end{equation*}
called the \textit{index} of the family. Notice that the rank map $K(X)\rightarrow\mathbb{Z}$ sends the index of the family to the index of a single operator in the family. It turns out that $\mathrm{ind}\{T_x\}$ is independent of choices, and that the assignment is a bijection. A key input in the proof is Kuiper's theorem, which states that the unitary group of a Hilbert space is contractible.\footnote{This also implies that any Hilbert bundle over a CW complex is trivial. In particular, after trivializing, a family of operators on a Hilbert bundle can be considered as a family of operators on fixed Hilbert space; we implicitly did so for the family of twisted Dirac operators on $S^1$.}
\begin{remark}\label{index0}
For our purposes, it is important to point out the following: if $\{T_x\}$ is a family of isomorphisms, the corresponding element in $K(X)$ is zero. This is because we can choose $W=\{0\}$.
\end{remark}

\vspace{0.3cm}

\subsection{The index theorem for families of twisted Dirac operators}
Consider now the family of twisted Dirac operators $\{D_B^+\}$ on $M$ parametrized by the torus $\mathbb{T}_M$, and the corresponding element
\begin{equation*}
\mathrm{ind}\{D_B^+\}\in K(\mathbb{T}_M).
\end{equation*}
Theorem \ref{AS} gives us a formula for the image of such element under the dimension map $K(\mathbb{T}_M)\rightarrow \mathbb{Z}$; we would like to obtain more detail information. Of course, $K$-theory is a quite complicated object to deal with, so it would be desirable to work with something simpler such as cohomology. Luckily, one can use Chern classes $c_i$ to define a ring homomorphism
\begin{equation*}
\mathrm{ch}: K(X)\rightarrow H^{\mathrm{even}}(X;\mathbb{Q})
\end{equation*}
called the Chern character; this is given on a vector bundle $E$ by
\begin{equation*}
\mathrm{ch}(E)=\mathrm{dim}E+c_1(E)+\frac{c_1^2-2c_2}{2}+\dots;
\end{equation*} 
in particular, the degree zero part of $\mathrm{ch}$ is exactly the dimension. It turns out that $\mathrm{ch}$ is a rational isomorphism. 
\begin{remark}
On the other hand, the torsion subgroups of $K(X)$ and $H^{\mathrm{even}}(X;\mathbb{Z})$ are in general very different.
\end{remark}
The index theorem for families will provide a formula for $\mathrm{ch}(\mathrm{ind}\{D_B^+\})\in H^{\mathrm{even}}(\mathbb{T}_M;\mathbb{Q})$ in terms of the topology of $M$. Recall that $\mathbb{T}_M=H^1(M;\mathbb{R})/H^1(M;2\pi\mathbb{Z})$, so that we can identify
\begin{equation*}
H_1(\mathbb{T}_M;\mathbb{Z})=H^1(M;\mathbb{Z}).
\end{equation*}
Via this identification, a basis $x_1,\dots ,x_m$ of $H^1(M;\mathbb{Z})$ determines a dual basis $y_1,\dots,y_m$ of $H^1(\mathbb{T}_M;\mathbb{Z})=H_1(\mathbb{T}_M;\mathbb{Z})^*$. We can then consider the element
\begin{equation*}
\Omega=\sum x_i\otimes y_i\in H^1(M;\mathbb{Z}) \otimes H^1(\mathbb{T}_M;\mathbb{Z})\subset H^2(M\times \mathbb{T}_M;\mathbb{Z}),
\end{equation*}
which is independent of the choice of basis\footnote{More intrinsically, $\Omega$ is the first Chern class of the so-called \textit{Poincar\'e line bundle} on $M\times \mathbb{T}_M$.}. Recall that the A-hat genus is the combination of Pontryagin classes
\begin{equation*}
\widehat{\mathcal{A}}(M)=1-\frac{p_1}{24}+\frac{7p_1^2-4p_2}{5760}+\dots,
\end{equation*}
which is an element of $H^{\mathrm{even}}(M;\mathbb{Q})$.
\begin{thm}[\cite{ASIV}]\label{indexfamily}
Consider an even dimensional spin manifold $M$. Then
\begin{equation*}
\mathrm{ch}(\mathrm{ind}\{D_B^+\})=\langle e^{\Omega}\cdot \widehat{\mathcal{A}}(M),[M]\rangle\in H^{\mathrm{even}}(\mathbb{T}_M.\mathbb{Q}).
\end{equation*}
Here $e^{\Omega}=\sum_{n\geq0} \Omega^n/n!$ and $e^{\Omega}\cdot \widehat{\mathcal{A}}(M)$ defines an element in $H^{\mathrm{even}}(M\times \mathbb{T}_M;\mathbb{Q})$, which can be evaluated on the fundamental class $[M]$ to obtain an element in $H^{\mathrm{even}}(\mathbb{T}_M;\mathbb{Q})$.
\end{thm}
We refer the reader to \cite{LM} for a more detailed discussion of this result and its generalizations.

\vspace{0.3cm}
\subsection{Metrics of positive scalar curvature on the torus}
We are now ready to prove the following theorem of Gromov and Lawson.
\begin{thm}[\cite{GL}]
The torus $T^n$ does not admit a metric of positive scalar curvature.\end{thm}
\begin{remark}
In fact, they prove the stronger result that a metric on $T^n$ with $s\geq0$ is necessarily flat.
\end{remark}
To see this, we can assume $n$ is even (otherwise we can consider $S^1\times T^n$). If $s>0$, then the argument of Proposition \ref{Bochner} shows that $D_B^+$ has trivial kernel and cokernel, hence is an isomorphism, for all $B$. By Remark \ref{index0}, this implies that $\mathrm{ch}(\mathrm{ind}\{D_B^+\})=0$. As $\widehat{\mathcal{A}}(T^n)=1$, the index formula simplifies to 
\begin{align*}
\mathrm{ch}(\mathrm{ind}\{D_B^+\})&=\langle e^{\Omega},[T^n]\rangle=\langle \frac{\Omega^n}{n!},[T^n]\rangle=\\
&=\pm\langle x_1\cup\cdots \cup x_n\otimes y_1\cup \cdots \cup y_n,[T^n]\rangle=\pm y_1\cup\cdots\cup y_n\neq 0
\end{align*}
because the cup products $x_1\cup\cdots\cup x_n$ and $ y_1\cup\cdots\cup y_n$ evaluate to $\pm1$ on the fundamental classes $[T^n]$ and $[\mathbb{T}_{T^n}]$ of the tori respectively, and we obtain a contradiction.

\vspace{0.3cm}

\subsection{Exercises}
\begin{enumerate}
\item Recall that $BU$ is simply connected and has $H^2(BU;\mathbb{Z})=\mathbb{Z}$, generated by the first Chern class $c_1$. Use this to compute $K(S^2)$, and describe explicitly the generators.
\item We will discuss a very concrete example of Kuiper's theorem. Consider the family of unitary maps $T_{\vartheta}$ on $l^2$ parametrized by $\vartheta\in S^1$ given in the standard basis $\{e_i\}$ by
\begin{equation*}
T(e_i)=\begin{cases}
e^{i\vartheta} e_i \text{ if }i=0\\
e_i\text{ if }i>0.
\end{cases}
\end{equation*}
Show that this this family is homotopically trivial by showing that the family
\begin{equation*}
\begin{bmatrix}
e^{i\vartheta} & 0\\
0 & e^{-i\vartheta}
\end{bmatrix}\in\mathrm{SU}(2)
\end{equation*}
is trivial, and applying the Mazur swindle
\begin{equation*}1=1+(-1+1)+(-1+1)+\dots=(1-1)+(1-1)+\dots=0.
\end{equation*}
\item Show that an orientable surface $\Sigma$ is spin. Use the index theorem for families to show that on any Riemannian surface $(M,g)$ of genus at least one there exists a twisted Dirac operator $D_B^+$ which has non-trivial kernel.
\end{enumerate}

\vspace{2cm}

\section{Dirac operators and the Weinstein conjecture in dimension three}

\subsection{Self-adjoint operators} In our previous examples we focused on the case of even dimensional manifolds, as in this case we have the decomposition of $D_B$ into two chiral operators $D^{\pm}_B$ adjoint to each other. In odd dimensions, there is no such canonical splitting, so we can only work with the Dirac operators $D_B$ themselves. Notice that as these are self-adjoint, $\mathrm{coker}D_B=\mathrm{ker}D_B$ so they all have automatically index zero, and hence there is no interesting analogue of Theorem \ref{AS}. On the other hand, it turns out that it is interesting to study families of them. More abstractly, consider again a complex separable Hilbert space $H$. Recall that an operator $T$ is \textit{self-adjoint} if 
\begin{equation*}
\langle Tx,y\rangle=\langle x,Ty\rangle\text{ for all }x,y\in H.
\end{equation*}
This is the infinite dimensional analogue of a being symmetric, but one has to be careful as the usual spectral theorem in general does not hold: in fact, $T$ might not even have eigenvalues at all. On the other hand, the differential operators we are interested in (e.g. $\Delta$ and $D_B$) are always diagonalizable on a compact manifold, so for simplicity we will only discuss their case\footnote{In general, one needs to rephrase the discussion that follows in terms of the quadratic form $x\mapsto \langle x,Tx\rangle$.}. It turns out that the spectral theory of such operators can have very different qualitative behavior; for example, looking at complex-valued functions on the circle $S^1$:
\begin{itemize}
\item for $c\in\mathbb{R}$ the shifted Laplacian $\Delta+c=-\frac{d^2}{dx^2}+c$ has eigenvalues $n^2+c$ for $n\in \mathbb{Z}$, hence only finitely many negative eigenvalues and infinitely many positive ones. We call operators like this \textit{essentially positive};
\item on the other hand $-\Delta+c$ has infinitely many negative eigenvalues and finitely many positive ones, and we call it \textit{essentially negative};
\item finally, the twisted Dirac operator $i\frac{d}{dx}+c$ has eigenvalues $n+c$ for $n\in \mathbb{Z}$ hence it has infinitely many negative eigenvalues and infinitely many positive eigenvalues. We call operators like this \textit{indefinite}.
\end{itemize}
We then have the following.

\begin{thm}[\cite{ASskew}] The space $\mathrm{Fred}^{\mathrm{sa}}(H)$ of Fredholm self-adjoint operators on $H$ has three connected components
\begin{equation*}
\mathrm{Fred}^{\mathrm{sa}}(H)=\{\text{essentially positive}\}\coprod\{\text{essentially negative}\}\coprod\{\text{indefinite}\}.
\end{equation*}
The first two components are contractible, while the last component has the homotopy type of $U=\bigcup_{n\geq1} U(n)$.
\end{thm}
\begin{remark}
One should not confuse $U$ with the unitary group $U(H)$ of a Hilbert space $H$; the latter is contractible by Kuiper's theorem, while for example $\pi_1(U)=\mathbb{Z}$.
\end{remark}
As in the case of the circle $S^1$, a twisted Dirac operator $D_B$ on a compact manifold $M$ is an indefinite operator, hence the family $\{D_B\}$ gives rise to a map
\begin{equation*}
\mathbb{T}_M\rightarrow U,
\end{equation*}
well-defined up to homotopy.
\par
Now, for a finite connected CW complex $X$, the set $[X,U]$ can be identified with $K^1(X)$ the \textit{odd} $K$-theory of $X$; we also denote $K(X)$ by $K^0(X)$, and call it even $K$-theory.
\begin{remark}
More concretely, $K^1(X)$ is usually defined as the the kernel of the rank homomorphism $K^0(\Sigma X)\rightarrow \mathbb{Z}$, where $\Sigma X$ is the reduced suspension of $X$. This is equivalent to the definition above because 
\begin{equation*}
\mathrm{ker}\left(K^0(\Sigma X)\rightarrow \mathbb{Z}\right)=[\Sigma X, BU]=[X ,\Omega BU]=[X,U],
\end{equation*}
where we used that reduced suspension $\Sigma$ is the adjoint of loopspace $\Omega$, and $\Omega BG$ is equivalent to $G$ for any $G$.
\end{remark}
In particular, a family $\{S_x\}$ of indefinite self-adjoint operators parametrized by $X$ determines an index element in $\mathrm{ind}\{S_x\}\in K^1(X)$. As in the even case, there is a Chern character map
\begin{equation*}
\mathrm{ch}:K^1(X)\rightarrow H^{\mathrm{odd}}(X;\mathbb{Q})
\end{equation*}
which is a rational isomorphism. We are then ready to state the analogue of Theorem \ref{indexfamily} in the case of odd dimensional manifolds.
\begin{thm}[\cite{APS3}]\label{indexodd}
Consider an odd dimensional spin manifold $M$. Then
\begin{equation*}
\mathrm{ch}(\mathrm{ind}\{D_B\})=\langle e^{\Omega}\cdot \widehat{\mathcal{A}}(M),[M]\rangle\in H^{\mathrm{odd}}(\mathbb{T}_M;\mathbb{Q})
\end{equation*}
where again $e^{\Omega}=\sum \Omega^n/n!$.
\end{thm}
For example, when $M=S^1$, this shows that the family of twisted Dirac operators $i\frac{d}{dt}+c$ corresponds to a generator of $\pi_1(U)=\mathbb{Z}$.
\begin{remark}
In fact,  $\pi_1(U)=\mathbb{Z}$ is very closely related to the \textit{spectral flow} of a family of operators, a key concept in Floer theory.
\end{remark}

\vspace{0.3cm}

\subsection{The Weinstein conjecture} We want to discuss how families of Dirac operators are related to a dynamical problem in contact geometry that we now recall (see \cite[Part IV]{Can} for a more detailed discussion). Recall that, given an oriented $(2n-1)$-dimensional manifold $Y$, a \textit{contact form} is a $1$-form $\alpha$ such that
\begin{equation*}
\alpha\wedge (d\alpha)^{n-1}>0 \text{ at each point.}
\end{equation*}
The basic example is $\mathbb{R}^{2n-1}$ with coordinates $(x_1,y_1,\dots, x_{n-1},y_{n-1},z)$ and $\alpha= dz-\sum y_idx_i$. Now, at each point $d\alpha$ is a skew-symmetric bilinear form with maximal rank, so there exists a unique vector field $R_{\alpha}$, called the \textit{Reeb vector field}, for which
\begin{equation*}
d\alpha(R_{\alpha},-)=0\text{  and  } \alpha(R_{\alpha})=1.
\end{equation*}
For example, in the case of $\mathbb{R}^n$, $R_{\alpha}=\frac{d}{dz}$. The following conjecture is motivated by questions in classical (Hamiltonian) mechanics.
\begin{conj}[Weinstein] If $(Y,\alpha)$ is a closed contact manifold, then the Reeb vector field $R_{\alpha}$ admits a closed orbit.
\end{conj}
Of course, the example of $\mathbb{R}^n$ above shows the assumption that $M$ is closed is necessary. The Weinstein conjecture has been proved in many cases, under assumptions with very different nature; we refer for \cite{Hut} for a more thorough discussion of the literature. Our goal here is to discuss how families of Dirac operators enter in the proof of the following celebrated result of Taubes.
\begin{thm}[\cite{Tau}] The Weinstein conjecture holds in dimension three.
\end{thm}
The rough idea is the following: one can perturb the (curvature part of the) Seiberg-Witten equations on $Y$ \cite[Ch. 4]{KM} by a very large multiple of the contact form $\alpha$; in this limit, solutions to the Seiberg-Witten equations naturally give rise to certain union of orbits of $R_{\alpha}$ (including the \textit{empty} orbit). Because the count of solutions is independent of the perturbation (in a suitable sense), one can use computations done purely in Seiberg-Witten theory to show that there are closed Reeb orbits.
\par
In fact, the following result (which is really is a consequence of the subsequently proved result that embedded contact homology coincides with monopole Floer homology, see \cite{Hut}) holds.
\begin{thm}[Taubes] If $(Y^3,\alpha)$ has no closed Reeb orbit, then for the monopole Floer homology\footnote{This is phrased in terms of cohomology in \cite{Hut}; the statement here follows from the universal coefficient theorem.} groups of $Y$ we have
\begin{equation*}
\HMf_*(Y,\spin)=
\begin{cases}
\mathbb{Z} \text{ if }\spin\text{ is the spin$^c$ structure corresponding to }\mathrm{ker}(\alpha)\\
0\text{ for all other spin$^c$ structures on $Y$.}
\end{cases}
\end{equation*}
Here, $\mathbb{Z}$ corresponds to the empty Reeb orbit.
\end{thm}

The key step is then the following non-triviality result of Kronheimer and Mrowka.
\begin{thm} [\cite{KM}, Chapter 35] \label{nonvan}For any $Y$, if $\spin$ is a torsion spin$^c$ structure, $\HMf_*(Y,\spin)$ is non-vanishing (and in fact has positive rank) in infinitely many gradings.
\end{thm}
This concludes the argument because every three-manifold admits a torsion spin$^c$ structure (for example one induced by a spin structure). Recall that for a torsion spin$^c$ structure $\HMf_*(Y,\spin)$ correponds to the Heegaard Floer homology group $HF^-(Y,\spin)$ (see \cite{CGH}, \cite{KLT} and subsequent papers).

\vspace{0.3cm}
\subsection{Relations with the triple cup product}
Theorem \ref{nonvan} is proved by showing that the bar version of monopole Floer homology $\HMb_*(Y,\spin)$ has positive rank in infinitely many gradings. Recall that using the $\mathbb{Z}[U]$-module structure, this can be thought as the localization
\begin{equation*}
\HMf_*(Y,\spin)\otimes_{\mathbb{Z}[U]}\mathbb{Z}[U,U^{-1}],
\end{equation*}
and for a torsion spin$^c$ structure it corresponds to $HF^{\infty}(Y,\spin)$ in Heegaard Floer homology.
\par
The bar version turns out to be a much simpler object to study. Indeed, Kronheimer and Mrowka provide in \cite{KM} a complete computation of $\HMb_*(Y,\spin;\mathbb{Q})$. Later, in \cite{Lid} Lidman provided a computation with $\mathbb{Z}/2$-coefficients using surgery techniques in Heegaard Floer homology. Finally, in recent joint work with Miller Eismeier \cite{LME}, we built on the work of Kronheimer and Mrowka to provide a full computation with $\mathbb{Z}$-coefficients using techniques from homotopical algebra. We will give the idea of proof of the following statement, which is an important step in the proof of Theorem \ref{nonvan}.
\begin{thm}[\cite{KM}, Ch. 35] For any $Y$, if $\spin$ is a torsion spin$^c$ structure, $\HMb_*(Y,\spin)$ only depends on the triple cup product map
\begin{align*}
\cup^3_Y&: \Lambda^3H^1(Y;\mathbb{Z})\rightarrow \mathbb{Z}\\
\alpha\wedge\beta\wedge\gamma&\mapsto \langle \alpha\cup\beta\cup\gamma,[Y]\rangle.
\end{align*}
\end{thm}
In particular, $\HMb_*(Y,\spin)$ only depends on the algebraic topology of $Y$, while $\HMf_*(Y,\spin)$ is a much subtler invariant of $Y$.
Intuitively, this is because while $\HMf_*(Y,\spin)$ involves both irreducible and reducible solutions to the Seiberg-Witten equations, $\HMb_*(Y,\spin)$ only involves the reducible ones. Because of this, it can be described purely in terms of the family of Dirac operators $\{D_B\}$ parametrized by the torus $\mathbb{T}_Y$.
\begin{remark}
This is a manifestation in Floer homology of the classical localization theorem in $S^1$-equivariant cohomology, see for example \cite{Tu}.
\end{remark}
More specifically, associated with $(Y,\spin)$ there is a family of Dirac operators $\{D_B\}$ parametrized by $\mathbb{T}_Y$ which generalizes the family of twisted Dirac operators associated to a spin structure.\footnote{Here one should think of $\mathbb{T}_Y$ as the space of reducible solutions to the Seiberg-Witten equations.} After choosing a Morse function $f:\mathbb{T}_Y\rightarrow \mathbb{R}$, the chain complex computing $\HMb_*(Y,\spin)$ is generated by pairs $(B_*, E)$ where $B_*$ is a critical point of $f$, and $E$ is an eigenspace of $D_{B_*}$. The differential counts solutions to equations combining flowlines $B(t)\in\mathbb{T}_Y$ of $-\nabla f$ and the family of operators $\{D_{B(t)}\}$.
\par
Then the key point is that this construction can be generalized to any family of indefinite self-adjoint Fredholm operators $\mathbb{T}_Y\rightarrow \mathrm{Fred}^{\mathrm{sa}}(H)$ (and not just the family of Dirac operators); the result, called \textit{coupled Morse homology}, only depends on the homotopy class of the family. In our case this is the index
\begin{equation*}
\mathrm{ind}\{D_B\}\in[\mathbb{T}_Y,U]=K^1(\mathbb{T}_Y).
\end{equation*}
For dimensional reasons we have $\widehat{A}(Y)=1$, so the index formula in Theorem \ref{indexodd} (which also holds in the slightly more general situation of torsion spin$^c$ structures) gives the computation
\begin{equation*}
\mathrm{ch}(\mathrm{ind}\{D_B\})=\cup^3_Y\in H^{\mathrm{odd}}(\mathbb{T}_Y)
\end{equation*}
under the identification
\begin{equation*}
H^3(\mathbb{T}_Y;\mathbb{Z})=\Lambda^3H^1(\mathbb{T}_Y;\mathbb{Z})=\left(\Lambda^3H^1(Y;\mathbb{Z})\right)^*.
\end{equation*}
Finally, one can show that $K^1(\mathbb{T}_Y)$ is torsion-free, so that the Chern character (which is a rational isomorphism) is injective and $\mathrm{ind}\{D_B\}$ is determined by its Chern character.
\\
\par
In fact, the same computation can be used to show that the classifying map $\mathbb{T}_Y\rightarrow U$ for the family $\{D_B\}$ factors through $\mathbb{T}_Y\rightarrow SU(2)$. This is a key ingredient in the proof of Theorem \ref{nonvan} because for these simpler families one can provide a concrete description of coupled Morse homology in terms of suitably `twisted' versions of more familiar theories. In particular, Theorem \ref{nonvan} is proved by identifying $\HMb^*(Y,\spin)$ over the reals with a version of \textit{twisted de Rham cohomology} for which the underlying vector space is $\Omega^*(\mathbb{T}_Y)\otimes\mathbb{R}[U,U^{-1}]$\footnote{Here we are working with cohomology, so the action of $U$ has degree $+2$.} with differential
\begin{equation*}
x\mapsto dx +(\zeta\wedge x)U^{-1}
\end{equation*}
for a closed $3$-form $\zeta$ with $[\zeta]=\cup^3_Y\in H^{3}(\mathbb{T}_Y;\mathbb{R})$ \cite[Chapter $34.3$]{KM}. Notice that this is a differential because $\zeta\wedge\zeta=0$.

\vspace{0.3cm}

\subsection{Exercises}

\begin{enumerate}
\item Let $f:S^1\rightarrow \mathbb{R}$ a non-constant continuous function. Show that the map on $L^2(S^1)$ (complex valued functions) given by multiplication by $f$ is self-adjoint but does not have eigenvalues.
\item Show that the family of operators $T_{\vartheta}$ in Exercise $2.2$ defines a loop in $U=\bigcup U(n)$, and in fact represents a generator of $\pi_1(U)$. Why does the Mazur swindle argument not apply in this situation?
\item Consider a link $L\subset S^3$. Show that the branched double cover $\Sigma(L)$ has vanishing triple cup product by considering the action of the covering involution on cohomology. Use this to show that the three-torus $T^3$ is not a branched double cover of $S^3$.
\end{enumerate}

\vspace{2cm}

\bibliographystyle{alpha}
\bibliography{biblio}

\end{document}